\documentclass[a4paper,12pt]{amsart}
\usepackage{amssymb}
\usepackage{ifthen}
\usepackage{graphicx}
\usepackage{color}
\nonstopmode
\newtheorem{thm}{Theorem}
\newtheorem{rem}{Remark}

\newtheorem{lem}{Lemma}

\newtheorem{claim}{Claim}

\newtheorem{conj}{Conjecture}
\newtheorem{prob}{Problem}

\theoremstyle{definition}
\newtheorem{defn}{Definition}

\newtheorem{case}{Case}

\newenvironment{pf}[1][]{%
 \vskip 1mm
 \noindent
 \ifthenelse{\equal{#1}{}}%
  {{\slshape Proof. }}%
  {{\slshape #1.} }%
 }%
{\qed\medskip}

\newcounter{alphabet}
\newcounter{tmp}
\newenvironment{Thm}[1][]{\refstepcounter{alphabet}%
\bigskip%
\noindent%
{\bf Theorem \Alph{alphabet}}%
\ifthenelse{\equal{#1}{}}{}{ (#1)}%
{\bf .} \itshape}{\vskip 8pt}

\makeatletter
\newcommand{\Ref}[1]{\@ifundefined{r@#1}{}{\setcounter{tmp}{\ref{#1}}\Alph{tmp}}}
\makeatother

\newenvironment{Prob}[1][]{\refstepcounter{alphabet}%
\bigskip%
\noindent%
{\bf Problem \Alph{alphabet}}%
{\bf .} \itshape}{\vskip 8pt}

\newcommand{\IR}{{\mathbb R}}

\newcommand{\IC}{{\mathbb C}}
\newcommand{\ID}{{\mathbb D}}

\newcommand{\real}{{\operatorname{Re}\,}}

\def\be{\begin{equation}}
\def\ee{\end{equation}}

\newcommand{\bee}{\begin{enumerate}}
\newcommand{\eee}{\end{enumerate}}

\newcommand{\blem}{\begin{lem}}
\newcommand{\elem}{\end{lem}}
\newcommand{\bthm}{\begin{thm}}
\newcommand{\ethm}{\end{thm}}

\newcommand{\bdefe}{\begin{defn}}
\newcommand{\edefe}{\end{defn}}
\newcommand{\bprob}{\begin{prob}}
\newcommand{\eprob}{\end{prob}}
\newcommand{\bques}{\begin{ques}}
\newcommand{\eques}{\end{ques}}
\newcommand{\bei}{\begin{itemize}}
\newcommand{\eei}{\end{itemize}}

\newcommand{\bde}{\begin{deter}}
\newcommand{\ede}{\end{deter}}
\newcommand{\bca}{\begin{case}}
\newcommand{\eca}{\end{case}}
\newcommand{\bcl}{\begin{claim}}
\newcommand{\ecl}{\end{claim}}
\newcommand{\bcon}{\begin{conj}}
\newcommand{\econ}{\end{conj}}
\newcommand{\bcons}{\begin{conjs}}
\newcommand{\econs}{\end{conjs}}
\newcommand{\bprop}{\begin{propo}}
\newcommand{\eprop}{\end{propo}}
\newcommand{\br}{\begin{rem}}
\newcommand{\er}{\end{rem}}
\newcommand{\brs}{\begin{rems}}
\newcommand{\ers}{\end{rems}}
\newcommand{\bo}{\begin{obser}}
\newcommand{\eo}{\end{obser}}
\newcommand{\bos}{\begin{obsers}}
\newcommand{\eos}{\end{obsers}}
\newcommand{\bpf}{\begin{pf}}
\newcommand{\epf}{\end{pf}}
\newcommand{\ba}{\begin{array}}
\newcommand{\ea}{\end{array}}
\newcommand{\beq}{\begin{eqnarray}}
\newcommand{\beqq}{\begin{eqnarray*}}
\newcommand{\eeq}{\end{eqnarray}}
\newcommand{\eeqq}{\end{eqnarray*}}

\newcounter{minutes}
\divide\time by 60
\newcounter{hours}
\multiply\time by 60 \addtocounter{minutes}{-\time}

\begin{document}
\title[Note on the convolution of harmonic mappings]
{Note on the convolution of harmonic mappings}

\thanks{
File:~\jobname .tex,
          printed: \number\day-\number\month-\number\year,
          \thehours.\ifnum\theminutes<10{0}\fi\theminutes}

\author{Liulan Li and Saminathan  Ponnusamy $^\dagger $
}
\address{Liulan Li, College of Mathematics and Statistics
 (Hunan Provincial Key Laboratory of Intelligent Information Processing and Application),
Hengyang Normal University, Hengyang,  Hunan 421002, People's
Republic of China} \email{lanlimail2012@sina.cn}

\address{S. Ponnusamy,
 Department of Mathematics,
Indian Institute of Technology Madras, Chennai-600 036, India.}
\email{samy@iitm.ac.in}



\thanks{The work of the first author is supported by NSF of China (No.
11571216), the Applied Characteristic Discipline Program in Hunan
Province, the Science and Technology Plan Project of Hunan Province
(No. 2016TP1020) and the Science and Technology Plan Project of
Hengyang City (2017KJ183). The  work of the second author is
supported by Mathematical Research Impact Centric Support of Department of Science and Technology,
India (MTR/2017/000367).
}
\maketitle

\begin{abstract}
Dorff et al. \cite{DN} formulated a question concerning the convolution of two right half-plane mappings,
where the normalization of the functions was considered incorrectly. In this paper,
we have reformulated the open problem in correct form and provided
a solution to it in a more general form. In addition, we also obtain two
new theorems which correct and improve some other results.
\end{abstract}

\maketitle \pagestyle{myheadings}
\markboth{Liulan Li and Saminathan  Ponnusamy}{Convolutions of harmonic mappings convex in one direction}

\section{Introduction}

In the recent years, the topic of harmonic mappings is one of the
most studied subject in geometric function theory of one and several
complex variables (see \cite{ABR-book,BL,Du}). In \cite{AC-12},
Aleman and Constantin used harmonic mappings to provide a new
approach towards obtaining explicit solutions to the incompressible
two-dimensional Euler equations. Moreover, while the general
solution is not available in explicit form, structural properties of
the system permit them to identify several families of explicit
solutions. More recently, Constantin and  Martin \cite{CM-17}
investigated the explicit solutions of two-dimensional
incompressible Euler equations, by using the Lagrangian coordinates,
which results in improving the work of \cite{AC-12}, and it is
interesting to observe that they are  related to harmonic mappings.

In this article, we will consider complex-valued harmonic mappings
$f$ defined on the open unit disk ${\mathbb D}=\{z \in {\mathbb
C}:\, |z|<1\}$, which have the canonical representation of the form
$f=h+\overline{g}$, where $h$ and $g$ are analytic in $\ID$. This
representation is unique with the condition $g(0)=0$. In terms of
the canonical decomposition of $f$, the Jacobian $J_f$ of
$f=h+\overline{g}$ is given by $J_f(z)=|h'(z)|^2 -|g'(z)|^2$.
According to the Inverse Mapping Theorem, if the Jacobian of a $C^1$
mapping from $\ID$ to $\IC$ is different from zero, then the
function is locally univalent. The classical result of Lewy implies
that the converse of this statement also holds for harmonic
mappings.  Thus, every harmonic function $f$ on $\ID$ is locally
one-to-one and sense-preserving on $\ID$ if and only if $J_f(z)>0$
in $\ID$, i.e. $|h'(z)|>|g'(z)|$ in $\ID$. The condition $J_f(z)> 0$
is equivalent to the existence of an analytic function $\omega _f $
in $\ID$ such that
$$ |\omega _f (z)|<1 ~\mbox{ for }~z\in \ID,
$$
where $\omega _f(z)=g'(z)/h'(z)$ is called the dilatation of $f$. When there is no confusion,  it is often
convenient to use $\omega $ instead $\omega _f$. Let ${\mathcal H}=\{f=h+\overline{g}:\, ~\mbox{$h(0)=g(0)=0$ and $h'(0)=1$} \}$.
The class ${\mathcal H}_0$ consists of those functions $f\in {\mathcal H}$ with $g'(0)=0$.

The family of all sense-preserving univalent harmonic mappings in
${\mathcal H}$ will be denoted by ${\mathcal S}_H$, and let
$\mathcal{S}_{H}^{0}={\mathcal S}_H\cap {\mathcal H}_0$. Clearly,
the familiar class ${\mathcal S}$ of normalized analytic univalent
functions in $\ID$ is contained in ${\mathcal S}_H^{0}$. The class
${\mathcal S}_H$ together with its geometric subclasses have been
studied extensively  by Clunie and Sheil-Small
\cite{Clunie-Small-84} and investigated subsequently by several
others (see \cite{Du} and the survey article \cite{SaRa2013}). In
particular, we consider the convolution properties of the class
${\mathcal K}_H$ (resp.  ${\mathcal K}_H^0$) of functions ${\mathcal
S}_H$ (resp. ${\mathcal S}_H^0$) that map the unit disk $\ID$ onto a
convex domain. For different properties of functions in ${\mathcal
K}_H^0$, we refer the readers to \cite{Clunie-Small-84,Du,Good02,SaRa2013}
and also to the recent articles \cite{FHM,HM,LiuPo}.


\subsection{Preliminaries and convex harmonic mappings}
One of the important and interesting geometric subclass of
${\mathcal S}_H$ which attracted function theorists is the class of
univalent harmonic functions $f$ for which the range $D=f(\ID)$ is
convex in the direction $\alpha$ $(0\leq \alpha< \pi)$, meaning the
intersection of $D$ with each line parallel to the line through $0$
and $e^{i\alpha}$ is an interval or the empty set (see, for example,
\cite{Clunie-Small-84,Do,DN,HengSch70,HengSch73,M,RZ}). Convex in
the direction $\alpha =0$ (resp. $\alpha =\pi/2$) is referred to as
convex in the horizontal (resp. vertical) direction.

It is known \cite{Clunie-Small-84} that a harmonic mapping
$f=h+\overline{g}$ belongs to $\mathcal{K}_{H}^{0}:={\mathcal
K}_H\cap {\mathcal H}_0$ if and only if, for each $\alpha \in [0,
\pi )$, the function $F=h-e^{2i\alpha}g$ belongs to $\mathcal{S}$
and is convex in the direction $\alpha$.

\bdefe
A function  $f=h+\overline{g} \in {\mathcal S}_H$ is said to be a
slanted half-plane mapping with $\gamma$ ($0\leq\gamma<2\pi$) if $f$
maps $\ID$ onto  $H_\gamma :=\{w:\,{\rm Re\,}(e^{i\gamma}w)
>-(1+a)/2\}$, where $-1<a<1$.
\edefe

Using the shearing method due to
Clunie and Sheil-Small \cite{Clunie-Small-84}  and the Riemann
mapping theorem, it is easy to see that such a mapping has the form
\begin{equation}\label{li4-eq2}
h(z)+e^{-2i\gamma}g(z)=\frac{(1+a)z}{1-e^{i\gamma}z}.
\end{equation}
Note that $h(0)=g(0)=h'(0)-1=0$ and $g'(0)=e^{2i\gamma}a$. The class
of all \textit{slanted half-plane mappings with $\gamma$} is denoted by
${\mathcal S}(H_{\gamma})$, and we denote by $\mathcal{S}^{0}(H_{\gamma})$  the subclass of
$\mathcal{S}(H_{\gamma})$ with $a=0$.
Obviously, each $f\in {\mathcal S}(H_{\gamma})$ (resp. $\mathcal{S}^{0}(H_{\gamma})$)
belongs to the convex family $\mathcal{K}_{H}$ (resp. $\mathcal{K}_{H}^{0}$). Evidently there are infinitely many
\textit{slanted half-plane mappings} with a fixed $\gamma$. It is
worth recalling that functions $f\in{\mathcal S}(H_{\gamma})$ with
$\gamma=0$ are usually referred to as \textit{the right half-plane mappings}, especially when $a=0$.
For example, if $f_{0}=h_{0}+\overline{g_{0}}$, where
\be\label{li4-eq2b}
h_{0}(z)=\frac{z-\frac{1}{2}z^2}{(1-z)^2}= \frac{1}{2}\left (\frac{z}{1-z}+\frac{z}{(1-z)^2} \right )
\ee
and
\be\label{li4-eq2c}
g_{0}(z)=\frac{-\frac{1}{2}z^2}{(1-z)^2}=\frac{1}{2}\left (\frac{z}{1-z}-\frac{z}{(1-z)^2} \right ),
\ee
then
$$
h_{0}(z)+g_{0}(z)=\frac{z}{1-z}
$$
showing that $f_{0}=h_{0}+\overline{g_{0}}\in {\mathcal S}^{0}(H_{0})$ with the dilatation $\omega_{0}(z)=-z$.
The function $f_0$ plays the role of extremal for many extremal problems for the convex family $\mathcal{K}_{H}^{0}$.

\subsection{Convolution of harmonic mappings}
For two harmonic mappings $f=h+\overline{g}$ and
$F=H+\overline{G}$ in $\mathcal{H}$ with power series of the form
$$f(z)
=z+\sum_{n=2}^{\infty}a_{n}z^{n}+\sum_{n=1}^{\infty}\overline{b_{n}}\overline{z}^{n}
~\mbox{ and }~
F(z)
=z+\sum_{n=2}^{\infty}A_{n}z^{n}+\sum_{n=1}^{\infty}\overline{B_{n}}\overline{z}^{n},
$$
we define the harmonic convolution (or Hadamard product) as follows:
$$(f*F)(z)=(h*H)(z)+\overline{(g*G)(z)}=z+\sum_{n=2}^{\infty}a_{n}A_{n}z^{n}+\sum_{n=1}^{\infty}\overline{b_{n}B_{n}}\overline{z}^{n}.$$

Clearly, the space ${\mathcal H}$ is closed under the operation
$\ast$, i.e. ${\mathcal H} \ast {\mathcal H}\subset {\mathcal H}$.
In the case of conformal mappings, the literature about convolution
theory is exhaustive (see for example \cite{samy95,PoSi96,rs1} and
the references therein). Unfortunately, most of these results do not
necessarily carry over to the class of univalent harmonic mappings
in $\ID$ (see \cite{Do,DN,Good02,LiPo1,LiPo2}).

What is surprising is that even if $f,\ F\in {\mathcal K}_H$, the
convolution $f\ast F$ is not necessarily locally univalent in $\ID$.
In view of this observation and some other reasonings, not much
could be achieved on the convolution of harmonic univalent mappings.
However, some progress has been obtained in the recent years, see
\cite{Do,DN,Good02,LiPo1,LiPo2}.

\subsection{Reformulation of the problem of Dorff et al. in \cite{DN} and main results}
In 2012,  Dorff et al. in \cite{DN} proved the following result.

\begin{Thm}{\rm (\cite[Theorem 2]{DN})}\label{ThmA}
If $f_k\in {\mathcal S^0}(H_{\gamma_k})$, $k=1, 2$, and $f_1\ast
f_2$ is locally univalent in $\ID$, then $f_1\ast f_2$ is convex in
the direction $-(\gamma_1 +\gamma_2)$.
\end{Thm}

By a similar reasoning, we can generalize the result to the setting
${\mathcal S}(H_{\gamma})$.

\blem\label{lem1}
If $f_k\in {\mathcal S}(H_{\gamma_k})$, $k=1,2$, and $f_1\ast f_2$
is locally univalent in $\ID$, then $f_1\ast f_2$ is convex in the
direction $-(\gamma_1 +\gamma_2)$.
\elem
\bpf Proof of this lemma follows by adopting the method of proof of \cite[Theorem 2]{DN}. So we omit the details.
\epf

Moreover, Dorff et al. in \cite{DN} also considered the situations where
$f_1\ast f_2$ is locally univalent and sense-preserving, and obtained

\begin{Thm}{\rm (\cite[Theorem 4]{DN})}\label{ThmB}
Let $f=h+\overline{g}\in  {\mathcal S}^{0}(H_{0})$ with dilatation
$\omega(z)=\frac{z+a}{1+az}$, where $a\in (-1,1)$. Then $f_0\ast
f\in {\mathcal S}_H^0$ and is convex in the horizontal direction.
\end{Thm}

Observe that for functions $f\in  {\mathcal S}^{0}(H_{0})$, the corresponding dilatation $\omega $ must satisfy the
condition that $\omega (0)=0$ which forces $a=0$ in Theorem \Ref{ThmB} and thus, Theorem \Ref{ThmB} is meaningful
only when $\omega(z)=z$. It is also worth recalling that, in 2010, Bshouty and Lyzzaik \cite{BL} brought out a
collection of open problems and conjectures on planar harmonic mappings, proposed
by many colleagues throughout the past quarter of a century.  In particular,
Dorff et al. \cite[Problem 3.26(a)]{BL}  posed the following open question.

\begin{Prob}\label{ProbC}
Let $f = h + \overline{g}\in {\mathcal S}^{0}(H_{0})$ with
dilatation $\omega(z)=(z +a)/(1+\overline{a}z)$, $|a|<1.$ Determine
other values of $a \in \mathbb{D}$ for which the  result of Theorem
{\rm \Ref{ThmB}} holds.
\end{Prob}


In Theorem \Ref{ThmB} and Problem \Ref{ProbC}, $f\in {\mathcal
S^0}(H_0)$ holds only when $a=0$. In other words, the normalization
was not taken care properly. In view of this reasoning, it is
essential to reconsider the above problem  in the setting ${\mathcal
S}(H_{\gamma})$ and hence, one has to reformulate Problem
\Ref{ProbC} taking into account of the normalization condition.

\bprob\label{prob2}
Let $f=h+\overline{g}\in {\mathcal S}(H_{\gamma})$ such that
$$h(z)+e^{-2i\gamma}g(z)=\frac{(1+a)z}{1-e^{i\gamma}z} ~\mbox{ and }~\omega(z)=e^{2i\gamma}\frac{ze^{i\theta}+a}{1+aze^{i\theta}}.
$$
Determine the values of $a$ and $\theta$ such that $f_0\ast f\in {\mathcal S}_H^0$ is univalent in $\ID$.
\eprob

Without realizing the error in Theorem \Ref{ThmB} and Problem
\Ref{ProbC},  the present authors \cite{LiPo1,LiPo2} investigated
the convolution properties of $f_0=h_0+\overline{g_0}$ with slanted
half-plane mappings $f\in {\mathcal S}(H_\gamma)$ and obtained
\cite[Theorem 2.2]{LiPo1} and \cite[Theorem 1.3]{LiPo2}. Moreover,
in a recent article, Liu and Ponnusamy \cite{LiuPo}, obtained the
corrected version of Theorem \Ref{ThmB} which we may now recall it
here.

\begin{Thm}{\rm (\cite[Theorem 1]{LiuPo})}\label{ThmC}
Let  $f=h+\overline{g}\in{\mathcal S}(H_{0})$ with
$$h+g= \frac{(1+a)z}{1-z} ~\mbox{ and }~
\omega(z)=\frac{z+a}{1+a z},
$$
where $-1<a<1$, and $f_1 =h_{1}+\overline{g_{1}}\in{\mathcal S^0}(H_{0})$ with dilatation
$\omega_{1}(z)=e^{i\theta}z~(\theta\in\mathbb{R})$. Then $f_{1}*f$
is locally univalent and convex in the horizontal direction.
\end{Thm}

In this paper, we determine a family of values $a$ and $\theta$ such
that the condition in Problem \ref{prob2} is satisfied (see Theorem
\ref{uni1}, which corrects in particular the two results of
\cite{LiPo1,LiPo2}, namely, \cite[Theorem 2.2]{LiPo1} and
\cite[Theorem 1.3]{LiPo2}). We also state and prove two new results,
namely, Theorems \ref{uni2} and \ref{uni3} which correct and improve
some others. Motivation and statements of these results are
discussed in Section \ref{sec-main} and their proofs will be given
in Section \ref{sec2}.

\section{Main Results}\label{sec-main}

\begin{thm}\label{uni1}
Let $f=h+\overline{g}\in {\mathcal S}(H_{\gamma})$ with
$$h(z)+e^{-2i\gamma}g(z)=\frac{(1+a)z}{1-e^{i\gamma}z} ~\mbox{ and }~
\omega(z)=e^{2i\gamma}\frac{ze^{i\theta}+a}{1+aze^{i\theta}},
$$
where $\theta\in\IR$ and $a\in (-1,1)$. If one of the following conditions holds, then
$f_0\ast f\in {\mathcal S}_H^0$ and is convex in the direction of
$-\gamma$:
\bee
\item [{\rm (1)}] $\cos(\theta-\gamma)=-1$ and $-1/3\leq a<1$.

\item [{\rm (2)}] $-1<\cos(\theta-\gamma)\leq1$ and $a^2<\frac{1}{5-4\cos(\theta-\gamma)}$.
\eee
\end{thm}

The proof of Theorem \ref{uni1} will be presented in Section \ref{sec2}.

In order to state and prove our next two results, we consider the
class  ${\mathcal S^0}(\Omega_\beta)$  of functions $f\in {\mathcal
S}_H^0$ such that $f$ maps $\ID$ onto the asymmetric vertical strip
domains
$$\Omega_\beta=\Big \{w:\, \frac{\beta-\pi}{2\sin\beta}<{\rm Re\,}w<\frac{\beta}{2\sin\beta}\Big \},
$$
where $0<\beta<\pi$. Each $f=h+\overline{g} \in {\mathcal S^0}(\Omega_\beta)$ has the form
\be\label{li4-eq3}
h(z)+g(z)=\psi (z), \quad \psi(z)=
\frac{1}{2i\sin\beta}\log\left(\frac{1+ze^{i\beta}}{1+ze^{-i\beta}} \right).
\ee

In \cite{KDGS, KGSD1, KGSD2}, Kumar et al. considered mappings
$F_a=H_a+\overline{G_a}$ with \be\label{li4-eq12}
H_a(z)+G_a(z)=\frac{z}{1-z}  ~\mbox{ and }~
\frac{G'_a(z)}{H'_a(z)}=\frac{a-z}{1-az} \ee and obtained some
convolution results of such mappings with mappings in ${\mathcal
S}^{0}(H_{0})\bigcup{\mathcal S^0}(\Omega_\beta)$. Again, from
\eqref{li4-eq12}, we see that $H_a'(0)+G_a'(0)=1$, i.e. $G_a'(0)=0$
which contradicts the second condition in \eqref{li4-eq12} unless
$a=0$. Before we state and prove corrected versions of their
results, we recall them.

\begin{Thm}{\rm (\cite[Theorem 2.2]{KDGS})}\label{ThmD}
Let  $f=h+\overline{g}\in{\mathcal S}^{0}(H_{0})$ with $h+g=z/(1-z)$
and dilatation $\omega(z)=e^{i\theta}z^n$, where $\theta\in\IR$ and
$n$ are positive integers. If $a\in[\frac{n-2}{n+2},1)$, then $f\ast
F_a$ is convex in the horizontal direction.
\end{Thm}

\begin{Thm}{\rm (\cite[Theorem 2.4]{KGSD1})}\label{ThmE}
Let $f=h+\overline{g}\in{\mathcal S^0}(\Omega_\beta)$ with
$\beta=\frac{\pi}{2}$ and dilatation $\omega(z)=e^{i\theta}z^n$,
where $\theta\in\IR$ and $n$ are positive integers.  If
$a\in[\frac{n-2}{n+2},1)$, then $f\ast F_a$ is convex in the
horizontal direction.
\end{Thm}

\begin{Thm}{\rm (\cite[Theorems 2.3, 2.5 and 2,6]{KGSD2})}\label{ThmF}
Let  $f=h+\overline{g}\in{\mathcal S^0}(\Omega_\beta)$ with
dilatation $\omega(z)=e^{i\theta}z^n$, where $0<\beta<\pi$,
$\theta\in\IR$ and $n$ are positive integers. If
$a\in[\frac{n-2}{n+2},1)$ and $n\leq4$, then $f\ast F_a$ is convex
in the horizontal direction.
\end{Thm}

As pointed out above, these results hold if and only if $a=0$. In order to reformulate these results in correct
form, we may consider harmonic mappings $f^a_0=h^a_0+\overline{g^a_0}$ such that
$$ h^a_0(z)+g^a_0(z)=\frac{(1+a)z}{1-z} ~\mbox{ and }~
\frac{(g^a_0)'(z)}{(h^a_0)'(z)}=\frac{a-z}{1-az}.
$$
In fact, we now consider a more general form of the above $f^a_0$ in the following way:
Let $f^a_\gamma=h^a_\gamma +\overline{g^a_\gamma}\in {\mathcal
S}(H_\gamma)$ with the dilatation
$$\omega (z)=-e^{2i\gamma}\frac{e^{i\gamma}z-a}{1-ae^{i\gamma}z},
~\mbox{ and
}~h^a_\gamma(z)+e^{-2i\gamma}g^a_\gamma(z)=\frac{(1+a)z}{1-e^{i\gamma}z}.
$$
Then a computation gives slanted half-plane mappings with $\gamma$ as $f^a_\gamma=h^a_\gamma+\overline{g^a_\gamma}$, where
$$h^a_\gamma(z)=\frac{(1+a)I_\gamma(z)+(1-a)zI'_\gamma(z)}{2},$$
$$g^a_\gamma(z)=e^{2i\gamma}\cdot\frac{(1+a)I_\gamma(z)-(1-a)zI'_\gamma(z)}{2}\;\
\mbox{and}\;\ I_\gamma(z)=\frac{z}{1-e^{i\gamma}z}.
$$
Obviously, when $\gamma=0$, $f^a_\gamma$ coincides with $f^a_0$
which has been considered above.

For any $f=h+\overline{g}\in\mathcal{H}$, the above representation
for $f^a_\gamma$ quickly gives that
$$\left(h^a_\gamma\ast h\right)(z)=\frac{(1+a)e^{-i\gamma}h(e^{i\gamma}z)+(1-a)zh'(e^{i\gamma}z)}{2},
$$
and
$$\left(g^a_\gamma\ast
g\right)(z)=e^{2i\gamma}\cdot\frac{(1+a)e^{-i\gamma}g(e^{i\gamma}z)-(1-a)zg'(e^{i\gamma}z)}{2}.
$$
Then by a computation, we see that the dilatation $\widetilde{\omega}$ of $f^a_\gamma\ast f$ is given by
\be\label{li4-eq4}
\widetilde{\omega(z)}=e^{2i\gamma}\cdot\frac{2ag'(ze^{i\gamma})-(1-a)e^{i\gamma}z
g''(ze^{i\gamma})}{2h'(ze^{i\gamma})+(1-a)ze^{i\gamma}h''(ze^{i\gamma})}.
\ee
For such slanted half-plane mappings $f^a_\gamma$, we obtain the following convolution theorems.

\begin{thm}\label{uni2}
Let $f=h+\overline{g}\in {\mathcal S^0}(H_{\gamma_1})$ with
$$h(z)+e^{-2i\gamma_1}g(z)=\frac{z}{1-e^{i\gamma_1}z}
$$
and dilatation $\omega(z)=e^{i\theta}z^n$, where $n$ are positive
integers and $\theta\in\IR$. If $a\in[\frac{n-2}{n+2},1)$, then
$f\ast f^a_\gamma$ is convex in the direction $-(\gamma_1+\gamma)$.
\end{thm}

\begin{thm}\label{uni3}
Let $f=h+\overline{g}\in{\mathcal S^0}(\Omega_\beta)$ with
dilatation $\omega(z)=e^{i\theta}z^n$, where $0<\beta<\pi$,
$\theta\in\IR$ and $n$ are positive integers. If
$a\in[\frac{n-2}{n+2},1)$, then $f\ast f^a_\gamma$ is convex in the
direction $-\gamma$.
\end{thm}

\begin{rem}
\bee
\item [{\rm (1)}] Theorem \ref{uni2} is the corrected version of Theorem \Ref{ThmD}.

\item [{\rm (2)}]
Theorem \ref{uni3} is not only the corrected version of Theorems \Ref{ThmE} and \Ref{ThmF}, but also a generalization of the two
theorems.
\eee
\end{rem}

\section{The proofs of Theorems \ref{uni1}, \ref{uni2} and \ref{uni3}}\label{sec2}

\subsection{Two lemmas} For the proof of Theorem \ref{uni1}, we need to prove a couple of lemmas.

\blem\label{dila1}
Let $f=h+\overline{g}\in  {\mathcal S}(H_\gamma)$ with dilatation
$\omega (z)=g'(z)/h'(z)$ and \be\label{li4-eq2a}
h(z)+e^{-2i\gamma}g(z)=\frac{(1+a)z}{1-e^{i\gamma}z}. \ee Then the
dilatation $\widetilde{\omega}$ of $f_0\ast f$ is \be\label{li4-eq5}
\widetilde{\omega}(z)=-ze^{-i\gamma}\left
(\frac{\omega^2(z)+e^{2i\gamma}[\omega(z)-\frac{1}{2}z\omega'(z)]
+\frac{1}{2}e^{i\gamma}\omega'(z)}
{1+e^{-2i\gamma}[\omega(z)-\frac{1}{2}z\omega'(z)]+\frac{1}{2}e^{-i\gamma}z^2\omega'(z)}\right
). \ee
\elem
\begin{proof} Assume the hypothesis that $f=h+\overline{g}\in  {\mathcal S}(H_\gamma)$ with $\omega (z)=g'(z)/h'(z)$. Then
$$g'(z)=\omega (z)h'(z)~\mbox{ and }~g''(z)=\omega'(z) h'(z)+\omega (z)h''(z).
$$
Moreover, as $h$ and $g$ are related by the condition \eqref{li4-eq2a}, the first
equality above gives
\be\label{li4-eq6}
h'(z)=\frac{1+a}{(1+e^{-2i\gamma}\omega(z))(1-e^{i\gamma}z)^2}
\ee
and therefore,
\be\label{li4-eq7}
h''(z)=(1+a)\left [\frac{-(1-e^{i\gamma}z)e^{-2i\gamma}\omega'(z)
+2(1+e^{-2i\gamma}\omega(z))e^{i\gamma}}{(1+e^{-2i\gamma}\omega(z))^2(1-e^{i\gamma}z)^3}\right ].
\ee
From the representation of $h_0$ and $g_0$ given by \eqref{li4-eq2b} and \eqref{li4-eq2c}, we see that
$$(h_0\ast h)(z)=\frac{h(z)+zh'(z)}{2}~\mbox{ and }~ (g_0\ast g)(z)=\frac{g(z)-zg'(z)}{2}.
$$
Therefore, as $f_0\ast f= h_0\ast h +\overline{g_0\ast g}$, the dilatation $\widetilde{\omega}$ of $f_0\ast f$ is given by
\be\label{li4-eq8}
\widetilde{\omega}(z)=\frac{(g_0\ast
g)'(z)}{(h_0\ast h)'(z)} =-\frac{zg''(z)}{2h'(z)+zh''(z)}=
-\frac{z\omega'(z) h'(z)+\omega (z)zh''(z)}{2h'(z)+zh''(z)}. \ee In
view of \eqref{li4-eq6} and \eqref{li4-eq7}, after some computation
\eqref{li4-eq8} takes the desired form. \end{proof}

\blem\label{dila2}
Let $f=h+\overline{g}\in {\mathcal S}(H_{\gamma})$ with
$$h(z)+e^{-2i\gamma}g(z)=\frac{(1+a)z}{1-e^{i\gamma}z} ~\mbox { and }~ \omega(z)=e^{2i\gamma}\frac{ze^{i\theta}+a}{1+aze^{i\theta}},
$$
where $\theta\in\IR$. Then the dilatation $\widetilde{\omega}$ of
$f_0\ast f$ is given by
$$\widetilde{\omega}(z)
=-ze^{3i\gamma}e^{2i\theta}\cdot\frac{(z+A)(z+B)}{(1+\overline{A}z)(1+\overline{B}z)},
$$
where \be\label{li4-eq9}
t(z)=z^2+\frac{3a+1}{2}e^{-i\theta}z+ae^{-2i\theta}+\frac{1-a}{2}e^{-i\gamma}e^{-i\theta},
\ee and $-A, -B$ are the two roots of $t(z)=0$ $($\,$A$ and $B$ may be equal$)$.
\elem
\bpf From the definition of $\omega$, it follows that
$$\omega'(z)=e^{2i\gamma}e^{i\theta}\frac{1-a^2}{(1+aze^{i\theta})^2}
$$
and thus,  by a computation, the expression for
$\widetilde{\omega}(z)$ in \eqref{li4-eq5} takes the form
\beqq
\widetilde{\omega}(z)
&=&-ze^{3i\gamma}e^{2i\theta} \left (\frac{z^2+\frac{3a+1}{2}e^{-i\theta}z+ae^{-2i\theta}+\frac{1-a}{2}e^{-i\gamma}e^{-i\theta}}
{1+\frac{3a+1}{2}e^{i\theta}z+ae^{2i\theta}z^2+\frac{1-a}{2}e^{i\gamma}e^{i\theta}z^2}\right )\\
& =&-ze^{3i\gamma}e^{2i\theta}\frac{t(z)}{t^*(z)},
 \eeqq
where $t(z)$ is given by \eqref{li4-eq9} and
$$t^*(z)=1+\frac{3a+1}{2}e^{i\theta}z+ae^{2i\theta}z^2+\frac{1-a}{2}e^{i\gamma}e^{i\theta}z^2.
$$
Suppose that $-A$ and $-B$ are the two roots of $t(z)=0$. Again, by
a simple calculation, it can be easily seen that
$$t^*(z)=z^2\cdot\overline{t(1/{\overline{z}})}
=z^2\cdot\overline{(1/{\overline{z}}+A)(1/{\overline{z}}+B)}=(1+\overline{A}z)(1+\overline{B}z).
$$
Therefore, the dilatation $\widetilde{\omega}$ of $f_0\ast f$ has the desired form.
\epf

\subsection{The proof of Theorem \ref{uni1}}
By Lemma \ref{lem1}, it suffices to prove that $f_1\ast f_2$ is
locally univalent and sense-preserving. So we only need to show that
the dilatation $\widetilde{\omega}$ of $f_0\ast f$ satisfies that
$$|\widetilde{\omega}(z)|<1,\;\ z\in\ID.
$$

Using the notation of Lemma \ref{dila2}, we write
$$t(z)=z^2+\frac{3a+1}{2}e^{-i\theta}z+ae^{-2i\theta}+\frac{1-a}{2}e^{-i\gamma}e^{-i\theta}=(z+A)(z+B).
$$

\bca The case $\cos(\theta-\gamma)=-1$ and $-1/3\leq a<1$.\label{case1} \eca

For this case,
$$t(z)=z^2+\frac{3a+1}{2}e^{-i\theta}z+\frac{3a-1}{2}e^{-2i\theta}=(z+e^{-i\theta})\left(z+\frac{3a-1}{2}e^{-i\theta}\right),
$$
so that the roots of $t(z)$ are $-A$ and $-B$, where $A=e^{-i\theta}$ and
$B=\frac{3a-1}{2}e^{-i\theta}$. Moreover,
$$t^*(z)=(1+e^{i\theta}z)\left(1+\frac{3a-1}{2}e^{i\theta}z\right)
$$
and from Lemma \ref{dila2} we find that
$$\widetilde{\omega}(z)=-ze^{3i\gamma}e^{2i\theta} \left (\frac{z+e^{-i\theta}}{1+e^{i\theta}z}\right)
\left
(\frac{z+\frac{3a-1}{2}e^{-i\theta}}{1+\frac{3a-1}{2}e^{i\theta}z}\right)
$$
which clearly implies that $|\widetilde{\omega}(z)|<1$ for $z\in\ID,$ since $-1/3\leq a<1$.

\bca The case $-1<\cos(\theta-\gamma)\leq1$ and
$a^2<\frac{1}{5-4\cos(\theta-\gamma)}$. \label{case2}\eca

Now let $a_0=ae^{-2i\theta}+\frac{1-a}{2}e^{-i\gamma}e^{-i\theta}$,
$a_1=\frac{3a+1}{2}e^{-i\theta}$ and $a_2=1$. Then
$t(z)=a_0+a_1z+a_2z^2$. A calculation yields that
$$|a_2|^2-|a_0|^2=\frac{1-a}{4}\cdot\left[a\left(5-4\cos(\theta-\gamma)\right)+3\right]>0$$
for  $-1<\cos(\theta-\gamma)\leq1$ and
$a^2<\frac{1}{5-4\cos(\theta-\gamma)}$. Now consider
$$t_1(z)=:\frac{\overline{a_2}t(z)-a_0t^*(z)}{z}=\frac{1-a}{4}\cdot\left[a\left(5-4\cos(\theta-\gamma)\right)+3\right]\cdot (z-z_0),
$$
where
$$z_0=\frac{(3a+1)(e^{-i\gamma}-2e^{-i\theta})}{a\left(5-4\cos(\theta-\gamma)\right)+3}=:\frac{u(a)}{v(a)},
$$
with
$$u(a)=(3a+1)(e^{-i\gamma}-2e^{-i\theta}) ~\mbox{ and }~  v(a)=a\left(5-4\cos(\theta-\gamma)\right)+3.
$$
A tedious calculation and the assumption show that
$$|v(a)|^2-|u(a)|^2=4\left(1+\cos(\theta-\gamma)\right)\left[1-a^2\left(5-4\cos(\theta-\gamma)\right)\right]
$$
which is positive by the assumption of Case \ref{case2}.
By using Cohn's rule (See for instance, \cite{RS2002}), the conclusion for this case follows.
The proof is complete. \hfill $\Box$


\subsection{The proof of Theorem \ref{uni2}} As in the proof of the previous theorem, by Lemma \ref{lem1},
it suffices to show that the dilatation $\widetilde{\omega}$ of $f\ast f^a_\gamma$
given by \eqref{li4-eq4} satisfies that $|\widetilde{\omega}(z)|<1$ for $z\in\ID.$
By \eqref{li4-eq4}, we see that $\widetilde{\omega(z)}=S(e^{i\gamma}z)$, where
$$S(z)=e^{2i\gamma}\left (\frac{2ag'(z)-(1-a)zg''(z)}{2h'(z)+(1-a)zh''(z)}\right ).
$$
Consequently, to complete the proof, it is enough to prove that $|S(z)|<1$ for $z\in\ID.$
The assumption that
$$h(z)+e^{-2i\gamma_1}g(z)=\frac{z}{1-e^{i\gamma_1}z}
$$
and the dilatation $\omega(z)=e^{i\theta}z^n$ yield that
$$g'(z)=e^{i\theta}z^n h'(z) ~\mbox{ and }~ g''(z)=ne^{i\theta}z^{n-1} h'(z)+e^{i\theta}z^nh''(z)
$$
so that $S(z)$ defined above takes the form
\be\label{li4-eq10}
S(z)=e^{(2\gamma+\theta)i}
z^n\left (\frac{2a-(1-a)n-(1-a)u[h(z)]}{2+(1-a)u[h(z)]}\right ),
\ee
where
\be\label{li4-eq11}
u[h(z)]=:z\frac{h''(z)}{h'(z)}=2\frac{ze^{i\gamma_1}}{1-ze^{i\gamma_1}}-n\frac{e^{(\theta-2\gamma_1)i}z^n}{1+e^{(\theta-2\gamma_1)i}z^n}.
\ee
Let $X=\real{u[h(z)]}$. Then $X>-1-\frac{n}{2}$ and thus, $2+n+2X>0$ for all $z\in\ID$.

By \eqref{li4-eq10}, it suffices to prove that for all $z\in\ID$,
$$T(z)=\left[2+(1-a)X\right]^2-\left[2a-(1-a)n-(1-a)X\right]^2\geq0.$$
By simplification, $T(z)$ reduces to
$$T(z)=(1-a)\left[2(1+a)-(1-a)n\right]\left[2+n+2X\right].$$
Finally, because $1-a>0$ and $2+n+2X>0$, we conclude that
$T(z)\geq0$ if and only if $a\in[\frac{n-2}{n+2},1)$. The desired conclusion follows.
\hfill $\Box$

\subsection{The proof of Theorem \ref{uni3}}
Recall that $f\ast f^a_\gamma$ is convex in the direction $-\gamma$
provided  $f\ast f^a_\gamma$ is locally univalent and
sense-preserving. By similar reasoning as in the proof of Theorem
\ref{uni2}, we only need to prove that for all $z\in\ID$,
$$T(z)=(1-a)\left[2(1+a)-(1-a)n\right]\left[2+n+2X\right]\geq0,
$$
where $X=\real{u[h(z)]}$ and $u[h(z)]$ is defined again by
$$u[h(z)]=z\frac{h''(z)}{h'(z)}.$$

On the other hand, the assumption that
$$h(z)+g(z)=\frac{1}{2i\sin\beta}\log\left(\frac{1+ze^{i\beta}}{1+ze^{-i\beta}} \right)
$$
and the dilatation $\omega(z)=e^{i\theta}z^n$ yield that
$$h'(z)=\frac{1}{(1+\omega(z))(1+ze^{i\beta})(1+ze^{-i\beta})},
$$
so that \beqq
u[h(z)]&=&-\frac{2(z+\cos\beta)z}{(1+ze^{i\beta})(1+ze^{-i\beta})}-\frac{z\omega'(z)}{1+\omega(z)}\\
&=&-\frac{2(z+\cos\beta)z}{(1+ze^{i\beta})(1+ze^{-i\beta})}-n\frac{e^{i\theta}z^n}{1+e^{i\theta}z^n}\\
&=&-\frac{e^{-i\beta}z}{1+e^{-i\beta}z}-\frac{e^{i\beta}z}{1+e^{i\beta}z}-n\frac{e^{i\theta}z^n}{1+e^{i\theta}z^n}.
\eeqq Clearly, the last equation implies that $2X>-2  -n$ for
$z\in\ID$. Finally, because $1-a>0$ and $(2+n)+2X>0$, it follows
that $T(z)\geq0$ if and only if $a\in[\frac{n-2}{n+2},1)$. The
desired conclusion follows. \hfill $\Box$

%
%

\end{document}